\documentclass{gtmon_a}
\pdfoutput=1

\proceedingstitle{The Zieschang Gedenkschrift}
\conferencestart{5 September 2007}
\conferenceend{8 September 2007}
\conferencename{Conference in honour of Heiner Zieschang}
\conferencelocation{Toulouse, France}

\editor{Michel Boileau}
\givenname{Michel}
\surname{Boileau}

\editor{Martin Scharlemann}
\givenname{Martin}
\surname{Scharlemann}

\editor{Richard Weidmann}
\givenname{Richard}
\surname{Weidmann}

\title{Roots of torsion polynomials and dominations}

\author{Michel Boileau}
\givenname{Michel}
\surname{Boileau}
\address{Institut de Math\'ematiques de Toulouse, UMR 5219\\
Universit\'e Paul Sabatier\\\newline 
31062 Toulouse Cedex 9
\\France\vspace{3pt}\\\newline
D\'ept de math\\UQAM\\\newline
PO Box 8888\\Centre-ville, Montr\'eal, Qc H3C 3P8\\Canada\vspace{3pt}
\\\newline
LMAM\\Department of Mathematics\\
Peking University\\\newline
Beijing 100871\\China}
\email{boileau@picard.ups-tlse.fr}
\urladdr{}

\author{Steve Boyer}
\givenname{Steve}
\surname{Boyer}
\email{boyer@math.uqam.ca}
\urladdr{}

\author{Shicheng Wang}
\givenname{Shicheng}
\surname{Wang}
\email{wangsc@math.pku.edu.cn}
\urladdr{}

\dedicatory{Dedicated to the memory of Heiner Zieschang}

\volumenumber{14}
\issuenumber{}
\publicationyear{2008}
\papernumber{05}
\startpage{75}
\endpage{81}

\doi{}
\MR{}
\Zbl{}

\arxivreference{math.GT/0606703}


\received{29 June 2006}
\revised{}
\accepted{5 February 2007}
\published{29 April 2008}
\publishedonline{29 April 2008}
\proposed{}
\seconded{}
\corresponding{}
\version{}

\makeatletter
\def\cnewtheorem#1[#2]#3{\newtheorem{#1}{#3}[section]
\expandafter\let\csname c@#1\endcsname\c@thm}
\makeatother

\AtBeginDocument{\let\bar\wbar\let\tilde\wtilde\let\hat\what}
\def\tildeW{\smash{\wwtilde W}\!}

\newtheorem{thm}{Theorem}[section]
\cnewtheorem{cor}[thm]{Corollary}
\cnewtheorem{lemma}[thm]{Lemma}
\cnewtheorem{prop}[thm]{Proposition}
\cnewtheorem{conj}[thm]{Conjecture}
\cnewtheorem{claim}[thm]{Claim}

\theoremstyle{definition}
\cnewtheorem{exa}[thm]{Example}
\cnewtheorem{que}[thm]{Question}
\cnewtheorem{rem}[thm]{Remark}
\cnewtheorem{defin}[thm]{Definition}
\cnewtheorem{add}[thm]{Addendum}
\cnewtheorem{obs}[thm]{Observation}

\def\pf{\begin{proof}}

\begin{document}

\begin{asciiabstract} 
We show that the nonzero roots of the torsion polynomials associated to the infinite cyclic covers of a given compact, connected, orientable 3-manifold M are contained in a compact part of C* a priori determined by M. This result is applied to prove that when M is closed, it dominates at most finitely many Sol manifolds. 
\end{asciiabstract}

\begin{htmlabstract}
We show that the nonzero roots of the torsion polynomials associated
to the infinite cyclic covers of a given compact, connected,
orientable 3&ndash;manifold M are contained in a compact part of
<b>C</b><sup>*</sup> a priori determined by M. This result is applied to
prove that when M is closed, it dominates at most finitely many
Sol manifolds.
\end{htmlabstract}

\begin{abstract} 
We show that the nonzero roots of the torsion polynomials associated
to the infinite cyclic covers of a given compact, connected,
orientable $3$--manifold $M$ are contained in a compact part of
$\mathbb{C}^*$ a priori determined by $M$. This result is applied to
prove that when $M$ is closed, it dominates at most finitely many
$Sol$ manifolds.
\end{abstract}

\maketitle

\section{Introduction}

All manifolds are connected and orientable in this paper. All homology groups will have $\mathbb Q$--coefficients unless otherwise specified.

Suppose that $M$ and $N$ are compact $3$--manifolds. We say that $M$ \textit{dominates} $N$ if there is a nonzero degree map $f \co (M, \partial M) \to (N, \partial N)$.

To each epimorphism $\psi \co \pi_1(M) \to \mathbb Z$ of the fundamental group of a compact $3$--manifold one can associate a torsion polynomial $\smash{\Delta_{\smash{\psi}}^M}(t)$. Our first result shows that the absolute values of the nonzero roots of such polynomials are pinched between two constants depending only on $M$, even though $\pi_1(M)$ has infinitely many epimorphisms to $\mathbb Z$ when its first Betti number is greater than one. We combine this result with a classical argument due to Wall for nonzero degree maps to show that the same conclusion holds for any $3$--manifold dominated by $M$. As an application we prove that  a closed $3$--manifold $M$ dominates at most finitely many $Sol$ manifolds.

\subsubsection*{Acknowledgements} The first author was partially supported
by the Institut de Math\'e\-ma\-tiques de Toulouse, UMR CNRS 5219.  
The second author was partially supported by NSERC.  The third author was partially supported by NSFC.

\section{Roots of torsion polynomials}

Given a compact $3$--manifold and an epimorphism $\psi \co \pi_1(M) \to \mathbb Z$, let  
$\smash{\wwtilde{M}_\psi} \to M$ be the associated infinite cyclic cover. The action 
on $\smash{H_1(\wwtilde{M}_\psi)}$ of $t = \smash{(T_\psi)_*}$ induced by the generator $T_\psi$ of the deck transformation group corresponding to $1 \in \mathbb Z$ makes $\smash{H_1(\wwtilde{M}_\psi)}$ a finitely generated $\Gamma$--module, where $\Gamma = \mathbb Q[\pi_1(M) / \hbox{ker}(\psi)] \cong  \mathbb Q[t, t^{-1}]$. Since 
$\Gamma$ is a principal domain, $\smash{H_1(\wwtilde{M}_\psi)} \cong \Gamma^{k}\oplus_{i = 1}^n \Gamma/(p_i(t))$ where $0 \ne p_i(t) \in \Gamma$. The product $\smash{\Delta_{\smash{\psi}}^M}(t) = p_1(t)p_2(t) \ldots p_n(t)$, called the \textit{torsion polynomial} of $\psi$, represents the order of the $\Gamma$--torsion submodule $\hbox{Tor}(H_1(\wwtilde{M}_\psi))$ of $\smash{H_1(\wwtilde{M}_\psi)}$ and is well-defined up to multiplication by some unit $r t^i$ of $\Gamma$ ($i \in \mathbb Z$, $0 \ne r \in \mathbb Q$). In particular, the set of nonzero roots $\{t_0 \in \mathbb C^* : \smash{\Delta_{\smash{\psi}}^M}(t_0) = 0\}$ is independent of the choice of $\smash{\Delta_{\smash{\psi}}^M}(t)$. A straightforward calculation shows that $\smash{\Delta_{\smash{\psi}}^M}(t)$ coincides, up to units, with the characteristic polynomial of the automorphism of the $\mathbb Q$--vector space $\oplus_{i = 1}^n \Gamma/(p_i(t))$ corresponding to multiplication by $t$.    

\begin{thm} \label{boundtorsionM}
A compact, connected, orientable $3$--manifold $M$ determines a constant $c_M > 0$ with the following property: If $t_0 \in \mathbb C^*$ is a root of a torsion polynomial $\smash{\Delta_{\smash{\psi}}^M}(t)$ associated to an epimorphism $\psi \co \pi_1(M) \to \mathbb Z$, then $\unfrac{1}{c_M} \leq |t_0| \leq c_M$. 
\end{thm}

\pf Since $1/t_0$ is a root of $\smash{\Delta_{\smash{-\psi}}^M}(t)$, it suffices to prove the existence of a constant $c_M$ such that $|t_0| \leq c_M$. 

For a group $G$ and $\alpha = \sum \limits_{g \in G} r_g g \in \mathbb Q[G]$, we set $\| \alpha \| = \sum \limits_{g \in G} |r_g|$. 

Consider a finite presentation $\langle x_j : r_i \rangle$ of $\pi_1(M)$ and let $J = \smash{\big(\frac{\partial r_i}{\partial x_j}\big)}$ be the associated Jacobian matrix. Define $k(\langle x_j : r_i \rangle) = \sum_{i,j} \smash{\big\| \frac{\partial r_i}{\partial x_j} \big\|} \in \mathbb N$ and set
$$k_M = \hbox{min} \{k(\langle x_j : r_i \rangle) : \langle x_j : r_i \rangle \hbox{ presents } \pi_1(M) \}.$$
We assume that $\langle x_j : r_i \rangle$ has been chosen to realize $k_M$ and that the number $m$ of generators is minimal among such presentations. 

Fix an epimorphism $\psi \co \pi_1(M) \to \mathbb Z$ and let $\Psi$ be the composition $\mathbb Z[\pi_1(M)] \to \mathbb Q[\pi_1(M)/ \hbox{ker}(\psi)] = \mathbb Q[t, t^{-1}]$. Recall that $J^\Psi = \smash{\big(\Psi\big(\frac{\partial r_i}{\partial x_j}\big)\big)}$ presents the $\Gamma$--module $H_1(\wwtilde{M}_\psi) \oplus \Gamma$ (see Burde and Zieschang \cite[Section 9]{BZ}, for example). Set $q_{ij}(t) = \smash{\Psi\big(\frac{\partial r_i}{\partial x_j}\big)} \in \mathbb Q[t, t^{-1}]$ and observe that $\| q_{ij}(t) \| \leq \smash{\big\| \frac{\partial r_i}{\partial x_j} \big\|}$. Thus the following claim holds. 

\begin{claim} \label{bound}
$\sum_{i,j} \| q_{ij}(t) \| \leq k_M$. 
\qed
\end{claim}

If $r$ denotes the $\Gamma$--rank of $J^\Psi$, then $\smash{\Delta_{\smash{\psi}}^M}(t)$ is, up to units, the g.c.d. of the $r$--rowed minors of $J^\Psi$ (see Jacobson \cite[Theorem 3.9]{J}, for example). Thus it suffices to show that the absolute values of the roots of some nonzero $r$--rowed minor of $J^\Psi$ are bounded above by a constant depending only on $M$. To that end, fix such a minor $D(t) \in \mathbb Z[t, t^{-1}]$ which, without loss of generality,  we can assume is polynomial in $t$, and let $D_0(t)$ be the monic polynomial with the same roots. Since $r \leq m$, the expansion of $D(t)$ in terms of the $q_{ij}(t)$ shows that $m ! k_M^m$ is an upper bound for the sum of the absolute values of its coefficients (cf \fullref{bound}). It is evident that the same inequality holds for $D_0(t) = t^s + b_{s-1} t^{s-1} + \ldots + b_0$. If $|t| > R = 1+ \sum_i |b_i|$, then $|D_0(t)| > R^n -(\sum_i |b_i| R^i) \geq R^{n-1}(R - \sum_i |b_i|) > 0$ so that the roots of $D_0(t)$ lie in the ball of radius $1+ \sum_i |b_i|$ centred at zero. Thus the theorem holds with $c_M = 1 +  m ! k_M^m$. 
\end{proof} 

We generalize this result with our applications in mind. 

\begin{thm} \label{boundtorsionN}
For a compact, connected, orientable $3$--manifold $M$, there is a constant $c_M > 0$ with the following property: If $N$ is a compact $3$--manifold dominated by $M$ and $t_0 \in \mathbb C^*$ is a root of a torsion polynomial $\smash{\Delta_{\smash{\psi}}^N}(t)$ of an epimorphism  $\psi \co \pi_1(N) \to \mathbb Z$, then $\unfrac{1}{c_M} \leq |t_0| \leq c_M$
\end{thm}

\pf  Suppose that $f\co M \to N$ is a nonzero degree map and fix an epimorphism $\psi \co \pi_1(N) \to \mathbb Z$. Since $\hbox{deg}(f) \ne 0$, there is an integer $n \geq 1$ such that the image $(\psi \circ f_\#)(\pi_1(M)) = n \mathbb Z$. Denote by $\theta\co \pi_1(M) \to \mathbb Z$ the epimorphism $(\unfrac{1}{n})(\psi \circ f_\#)$ and by $\Delta_\theta^M(t)$ the associated torsion polynomial. The theorem is a simple consequence of \fullref{boundtorsionM} and the following claim. 

\begin{claim} \label{nthpowerroot}
If $t_0 \in \mathbb C^*$ is a root of $\smash{\Delta_{\smash{\psi}}^N}(t)$, then $t_0^n$ is a root of $\Delta_\theta^M(t)$. 
\end{claim}

\proof Let  $\mathbb Q[t, t^{-1}]_f$ be the $\mathbb Z[\pi_1(M)]$--module whose underlying group is $\mathbb Q[t, t^{-1}]$ and whose $\pi_1(M)$ action is that determined by the homomorphism 
$f_\#\co \pi_1(M) \to \pi_1(N)$. Thus for $x \in \pi_1(M)$ and $p(t) \in \mathbb Q[t, t^{-1}]$ we have $x \cdot p(t) = t^{(\psi \circ f_\#)(x)} p(t)$. When $n = 1$, this action coincides with that of $\mathbb Z[\pi_1(M)]$ on $\mathbb Q[\pi_1(M) / \hbox{ker}(\psi \circ \smash{f_\#})]$ and so $H_1(M; \mathbb Q[t, t^{-1}]_f) \cong \smash{H_1(\wwtilde{M}_\theta)}$, where the latter has the $\Gamma$--action described above. In particular, since $\hbox{deg}(f) \ne 0$, there is a $\Gamma$--module splitting 
$$H_1(\wwtilde{M}_\theta) = H_1(M;\mathbb Q[t, t^{-1}]_f)) \cong H_1(N; \mathbb Q[t, t^{-1}]) \oplus K 
= H_1(\tilde{N}_\psi) \oplus K $$ 
for some finitely generated $\Gamma$--submodule $K$ of $H_1(\wwtilde{M}_\theta)$ (see the proof of \cite[Lemma 2.1]{Wall}). Hence when $n=1$, $\hbox{Tor}(H_1(\tilde{N}_\psi))$ is a $\Gamma$--submodule of $\hbox{Tor}(H_1(\wwtilde{M}_\theta))$, and so its order $\smash{\Delta_{\smash{\psi}}^N}(t)$ divides $\Delta_\theta^M(t)$, which implies the claim in this case. 

Next suppose $n > 1$ and let $\tilde{N}_{(\psi, n)} \to N$ be the $n$--fold cyclic cover with $\pi_1(\tilde{N}_{(\psi, n)})$ the kernel of the (mod $n$) reduction of $\psi$. Then $f$ lifts to a $\pi_1$--surjective, nonzero degree map $\tilde f\co M \to \tilde{N}_{(\psi, n)} = N'$.  Let $$\psi'\co\pi_1(N') \to n \mathbb Z \stackrel{\unfrac{1}{n}}{\longrightarrow}  \mathbb Z$$ be the epimorphism induced by $\psi$. The case $n = 1$ shows that any nonzero root of $\smash{\Delta_{\psi'}^{N'}(t)}$ is also a root of $\smash{\Delta_\theta^M(t)}$. On the other hand, it is easy to see that $\smash{(T_{\psi'})_* = (T_\psi)_{*}^n}$ on $\smash{H_1(\wwtilde{N'}_{\!\psi'})} = \smash{H_1(\tilde{N}_\psi)}$ so that if $t_0 \in \mathbb C^*$ is a root of $\smash{\Delta_{\smash{\psi}}^N}(t)$, then $\smash{t_0^n \in \mathbb C^*}$ is a root of $\Delta_{\psi'}^{N'}(t)$, and therefore of $\Delta_\theta^M(t)$. This completes the proof of the claim and therefore of \fullref{boundtorsionN}. 
\end{proof} 
 
\section[Q-Homology surface bundles]{$\mathbb Q$--Homology surface bundles}

Let $F$ be a compact surface and $A$ an abelian group. An \textit{$A$--homology} $F \times I$ is a $3$--manifold $W$ with boundary containing two disjoint surfaces $F_1 \cong F_2 \cong F$ such that
\begin{enumerate}
\item[(i)] $\overline{\partial W \setminus (F_1 \cup F_2)} \cong \partial F \times I$ where $\partial F \times \{0\} = \partial F_1, \partial F \times \{1\} = \partial F_2$, and 

\item[(ii)] the inclusion induced homomorphism $H_*(F_1; A) \to H_*(W; A)$ is an isomorphism. 
\end{enumerate}

(Duality and universal coefficients shows that (ii) is equivalent to each of the following three conditions: $H_*(W, F_1;A) = 0$; $H_*(W, F_2;A) = 0$; $H_*(F_2; A) \stackrel{\cong\ }{\smash{\longrightarrow}\vphantom{-}} H_*(W; A)$.) Note that $W$ determines orientations on $F_1$ and $F_2$ well-defined up to simultaneous  reversal. Thus the set $\hbox{Homeo}(F_2, F_1)^-$ of orientation reversing homeomorphisms $F_2 \to F_1$ is well-defined. For each $\varphi \in \hbox{Homeo}(F_2, F_1)^-$ we define $W_\varphi$ to be the compact, orientable manifold obtained from $W$ by identifying $F_2$ to $F_1$ via $\varphi$. The composition $$H_1(F_1;A) \xrightarrow{\cong} H_1(W;A) \xrightarrow{\cong}  H_1(F_2;A) \xrightarrow{\varphi_*} H_1(F_1;A)$$ determines an isomorphism 
$$\varphi_*^W\co H_1(F_1;A) \to H_1(F_1;A)$$ 
which we call the \textit{algebraic monodromy} of $W_\varphi$. Set 
$$\Delta_\varphi^W(t) = \det(\varphi_*^W - tI).$$
We call $W_\varphi$ an \textit{$A$--homology $F$ bundle}.

\begin{thm} \label{boundcharpoly}
For a compact, connected, orientable $3$--manifold $M$, there is a constant $c_M>0$ with the following property: If $W_\varphi$ is a $\mathbb Q$--homology surface bundle which is dominated by $M$, then the absolute values of the roots of the characteristic polynomial $\Delta_\varphi^W(t)$ of $\varphi_*^W$ are pinched between $1/c_M$ and $c_M$.
\end{thm}

\pf Let $F \subset W_\varphi$ be the nonseparating surface corresponding to $F_1 = \varphi(F_2)$. It determines a nonzero class $[F] \in H_2(W_\varphi)$, well-defined up to sign, and an epimorphism
$$\psi\co \pi_1(W_\varphi; \mathbb Z) \to \mathbb Z, \alpha \mapsto \alpha \cdot [F].$$ 
Let $\tildeW_\varphi \to W_\varphi$ be the infinite cyclic cover associated to this epimorphism $\psi$. 
Note that $H_1(\tildeW_\varphi) = H_1(W_\varphi; \Gamma)$ where $\Gamma$ is the $\mathbb Z[\pi_1(W_\varphi)]$--module $\mathbb Q[\pi_1(W_\varphi) / \hbox{ker}(\psi)] \cong  \mathbb Q[\mathbb Z] \cong \mathbb Q[t, t^{-1}]$. The $\mathbb Z[\pi_1(W_\varphi)]$ action on $H_1(\tildeW_\varphi)$ factors through one of $\Gamma$ in such a way that $t = (T_\varphi)_*$ where $T_\varphi\co \tildeW_\varphi \to \tildeW_\varphi$ is a  generator  of the group of deck transformations of $\tildeW_\varphi \to W_\varphi$. 

\begin{claim} \label{torsion} 
$H_1(\tildeW_\varphi)$ is a torsion module over $\Gamma$ whose order is represented by $\Delta_\varphi^W(t)$. 
\end{claim}

\pf  The quotient map $W \to W_\varphi$ lifts to an inclusion of $W$ into $\tildeW_\varphi$ with image $\tildeW_0$ say. Let $\smash{\tilde{F}_0} \subset \partial \tildeW_0$ correspond to $F_1$ and set $\tildeW_j = \smash{T_\varphi^j}(\tildeW_0), \smash{\tilde{F}_j} =  \smash{T_\varphi^j(\tilde{F}_0)}$. Then $\tildeW_\varphi = \cup_j \tildeW_j$ where $\tildeW_j \cap \tildeW_{k} = \emptyset$ if $|j-k| > 1$ and $\tildeW_j \cap \tildeW_{j-1} = \smash{\tilde{F}_j}$. Since $W$ is a $\mathbb Q$--homology $F_1 \times I$, the composition $H_1(F_1) = \smash{H_1(\tilde{F}_0)} \to H_1(\tildeW_\varphi)$ is an isomorphism under which the algebraic monodromy $\varphi_*^W\co H_1(F_1) \to H_1(F_1)$ corresponds to $(T_{\varphi})_*\co H_1(\tildeW_\varphi) \to H_1(\tildeW_\varphi)$. 

It is now clear that $H_1(\tildeW_\varphi)$ is a torsion module over $\Gamma$ since $H_1(\tildeW_\varphi) \cong H_1(F_1)$ is finite dimensional over $\mathbb Q$. Hence  the order of $H_1(\tildeW_\varphi)$ as a $\Gamma$--module corresponds to the characteristic polynomial of the automorphism $(T_{\varphi})_*$ of the $\mathbb Q$--vector space $H_1(\tildeW_\varphi)$, at least up to multiplication by some unit $\Gamma$. Since $(T_\varphi)_*$ corresponds to  $\varphi_*^W$ under $H_1(F_1) \stackrel{\cong\ }{\vphantom{-}\smash{\longrightarrow}} H_1(\tildeW_\varphi)$,  $\Delta_\varphi^W(t)$ also represents the order of $H_1(\tildeW_\varphi)$. 
\end{proof}

\fullref{torsion} shows that, up to multiplication by a unit, $\Delta_\varphi^W(t)$ is the torsion polynomial of the epimorphism $\psi$. \fullref{boundcharpoly} now follows from \fullref{boundtorsionN}. 
\end{proof}

\begin{cor} \label{finite} 
Let $W$ be a $\mathbb Q$--homology $F \times I$. A compact, connected, orientable $3$--manifold $M$ determines a finite subset ${\cal P}_{(M,W)}$ of $\mathbb Q[t]$ such that if $M$ dominates $W_{\varphi}$, then the characteristic polynomial of $\varphi_*^W$ is contained in ${\cal P}_{(M,W)}$. 
\end{cor}

\pf Let $\beta_1(F)$ be the first Betti number of $F$. The reader will verify that since $W$ is a $\mathbb Q$--homology $F \times I$, we can choose bases of for $H_1(F_1; \mathbb Z)$ and $H_1(F_2; \mathbb Z)$ with respect to which the matrix $X$ of $H_1(F_1) \to H_1(W) \to H_1(F_2)$ lies in $SL_{\beta_1(F)}(\mathbb Q)$ and the matrix $Y$ of $H_1(F_2) \stackrel{\varphi_*}{\longrightarrow} H_1(F_1)$ lies in $SL_{\beta_1(F)}(\mathbb Z)$. Now $\varphi^W_*$ is represented by $YX$ so the denominators of its entries are bounded above by some constant $N$. Thus the coefficients of $\smash{\Delta_\varphi^W(t)} = \det(YX - tI)$ have denominators bounded above by $\smash{N^{\beta_1(F)}}$ and since its degree is $\beta_1(F)$, the corollary follows from \fullref{boundcharpoly}. 
\end{proof}

\begin{rem} 
(1)\qua The finite set ${\cal P}_{(M,W)}$ described in the corollary depends on both $M$ and $W$. In the case when $W\cong F \times I$, the matrix $X$ of the proof of \fullref{finite} lies in $\smash{SL_{\beta_1(F)}}(\Z)$, so it is easy to see that ${\cal P}_{(M,W)}$ depends only on $M$ and the Euler characteristic of the fibre. 

(2)\qua A given compact $3$--manifold $M$ can be the total space of infinitely many distinct surface bundles over the circle. Moreover, there are cases where the Euler characteristic of the fibres are unbounded. However, \fullref{boundcharpoly} provides the following constraint on the monodromy of any such bundle structure on $M$.
\end{rem}

\begin{cor}
Given a compact $3$--manifold $M$, there is a constant $c_M > 0$ such that the absolute values of the roots of the characteristic polynomial of the algebraic monodromy of any surface bundle structure on $M\!$ are pinched between $\unfrac{1}{c_M}\!$ and $c_M$.~\qedsymbol
\end{cor}

Recall that an element $\varphi \in SL_2(\mathbb Z)$ is called \textit{hyperbolic} if $|\hbox{trace}(\varphi)| > 2$. 

\begin{cor} \label{finite1} 
A closed, connected, orientable $3$--manifold $M$ dominates only finitely many $Sol$ manifolds. 
\end{cor}

\pf First suppose that $M$ dominates a torus bundle over the circle with hyperbolic monodromy $\varphi \in SL_2(\mathbb Z)$. \fullref{finite} shows that there are only finitely many possibilities for $\hbox{trace}(\varphi)$, which is the negative of the coefficient of $t$ in $\Delta_\varphi^{T^2 \times I}(t)$. On the other hand, there are only finitely many $SL_2(\mathbb Z)$ conjugacy classes of hyperbolic elements of $SL_2(\mathbb Z)$ with a given trace (eg see Wang and Zhou \cite[Lemma 8]{WZ}). Since the homeomorphism type of a torus bundle over the circle depends only on the conjugacy class of its monodromy $\varphi \in SL_2(\mathbb Z)$, it follows that a closed, connected, orientable $3$--manifold can dominate at most finitely many torus bundles over the circle with hyperbolic monodromy. But a closed, connected $Sol$ manifold $N$ is double covered by such a bundle $\tilde N$ and so if $M$ dominates $N$, some double cover of $\wwtilde{M}$ dominates $\tilde N$. Since $M$ has only finitely many double covers, there are only finitely many possibilities for $\tilde N$, and therefore for $N$ \cite{Sa,MS}. 
\end{proof} 

It is known that if a closed, orientable $3$--manifold dominates a manifold which admits a geometric structure based on the geometries $\mathbb S^3, Nil$, or $\smash{\widetilde{SL_2}}$, then it dominates infinitely many distinct such manifolds \cite{W}. This is false for the remaining geometries. 

\begin{cor}
A closed, orientable $3$--manifold dominates at most finitely many manifolds admitting an $\mathbb S^2 \times \mathbb R, \mathbb E^3, \mathbb H^3, \mathbb H^2 \times \mathbb R$, or $Sol$ structure. 
\end{cor}

\pf The corollary holds for dominations of $\mathbb S^2 \times \mathbb R$ and $\mathbb E^3$ manifolds since there are only finitely many such spaces (see eg Scott \cite{Sc}). It holds for dominations of hyperbolic manifolds by Soma \cite{So}, for $\mathbb H^2 \times \mathbb R$ manifolds by 
Wang and Zhou \cite{WZ}, and for $Sol$ manifolds by \fullref{finite1}. 
\end{proof}

\bibliographystyle{gtart}
\bibliography{link}

\end{document}